\documentclass[5p,authoryear]{elsarticle}

\pdfoutput=1
\journal{}
\makeatletter
\def\ps@pprintTitle{%
   \let\@oddhead\@empty
   \let\@evenhead\@empty
   \let\@oddfoot\@empty
   \let\@evenfoot\@oddfoot
}
\makeatother

\usepackage{fancyhdr}

\usepackage[fontsize=9.2pt]{fontsize}
\usepackage{geometry}
\usepackage{natbib}
\usepackage[latin1]{inputenc}
\usepackage{amsmath}
\usepackage{amsfonts}
\usepackage{amssymb}
\usepackage{float} 
\usepackage{caption}
\usepackage{subcaption} 

\usepackage{amsthm} 
    \newtheorem{lemma}{Lemma}
    \newtheorem{assumption}{Assumption}
    \newtheorem{definition}{Definition}
    
    \newtheorem{remark}{Remark}
    \newtheorem{theorem}{Theorem}
    \newtheorem{corollary}{Corollary}
    \newtheorem{example}{Example}

\renewcommand\qedsymbol{$\blacksquare$} 

\usepackage[ruled, lined, ,longend, linesnumbered]{algorithm2e}
\usepackage[hang,flushmargin]{footmisc}
\usepackage{algpseudocode}
\usepackage{enumitem}
\SetKwInOut{Input}{Input}
\SetKwInOut{Declares}{Declares}
\SetKwInOut{Require}{Require}
\SetKw{KwBy}{by}
\SetKw{KwAnd}{and}
\SetKw{KwOr}{or}
\SetKw{KwEndFor}{end for}
\SetKwInput{KwParameters}{Parameters}
\SetKwInput{KwSelect}{Select}
\SetKwInput{KwNonDeterministic}{Non-deterministic}
\SetKwInOut{KwMySentence}{For each $k = 0, 1, 2, \dots$ repeat}
\SetKwRepeat{MyRepeat}{For each $k = 0, 1, 2, \dots$ repeat}{until}
\SetKwFor{MyRepeatEver}{For each $k = 0, 1, 2, \dots$ repeat}{}{}

\newlist{Rlist}{enumerate}{1}
\setlist[Rlist, 1]{%
    label=\textit{(\roman*)},
}

\newlist{RlistS}{enumerate}{1}
\setlist[RlistS, 1]{%
    label=\textit{(\roman*)},
    noitemsep,
}

\newlist{RlistL}{enumerate}{1}
\setlist[RlistL, 1]{%
    label=\textit{(\roman*)},
    labelindent=0pt,
    itemindent=1em,
    leftmargin=1em,
}

\newcommand{\R}{\mathbb{R}} 
\newcommand{\Rpp}{\R^{+}} 
\newcommand{\Z}{\mathbb{Z}} 
\newcommand{\T}{^\top} 
\newcommand{\cc}[1]{\mathcal{#1}} 
\def\sp#1#2{\langle #1,#2\rangle } 
\def\set#1#2{\{#1:#2\}} 
\def\fracg#1#2{{\displaystyle{\frac{#1}{#2}}}} 
\newcommand{\becauseof}[2][=]{\stackrel{\scriptstyle\mkern-1.5mu#2\mkern-1.5mu}{#1}} 
\def\Sum#1#2{\sum\limits_{#1}^{#2}} 
\def\tmin{{\rm min}} 
\def\tmax{{\rm max}} 
\def\st{{\rm s.t.}} 
\def\Ind{\cc{I}} 
\newcommand{\B}[2][]{%
\def\FirstArg{#1}
  \ifx\FirstArg\empty
      \cc{B}^n_{#2}
  \else
  \cc{B}^{#1}_{#2}
  \fi
} 

\newcommand{\fp}[1]{\hat{#1}}
\newcommand{\pq}[1][p.q]{_{{(#1)}}}
\def\fpx{\fp{x}}
\newcommand{\err}[1]{\texttt{err}\{#1\}} 
\newcommand{\exact}[1]{\texttt{exact}\{#1\}} 

\newcommand{\pP}{\cc{P}} 
\newcommand{\cQ}{\mathbb{Q}} 
\newcommand{\cP}{\mathbb{P}} 
\newcommand{\cX}{\cc{X}} 
\def\L{L} 
\def\ss{\sigma} 

\def\t{\rho} 
\def\ti{\t^{-1}} 
\newcommand{\CGM}[1][\t]{\cc{T}_{#1}} 
\def\Xinit{\cX^0} 
\def\fpXinit{\fp{\cX}^0} 
\newcommand{\g}[1]{\nabla f(#1)} 
\newcommand{\gfp}[1][\fpx^k]{\hat{g}_\t(#1)} 
\newcommand{\distfp}[1][\fpx^k]{\hat{d}^2(#1)} 

\def\tx{\tilde{x}}
\def\tu{\tilde{u}}

\def\stepk{s^k} 
\newcommand{\dopt}[1][k]{\gamma^{#1}}  
\makeatletter
\newcommand{\leqnomode}{\tagsleft@true}
\newcommand{\reqnomode}{\tagsleft@false}
\makeatother

\begin{document}
\pagestyle{fancy}
\thispagestyle{fancy}
\begin{frontmatter}

\title{\LARGE \bf Certification of the proximal gradient method under fixed-point arithmetic for box-constrained QP problems%
\tnoteref{t1}}
\tnotetext[t1]{This work was supported in part by Grant Margarita Salas, funded by the Ministerio de Universidades and the European Union (NextGenerationEU), by PRO3 MUR project Software Quality, and by PNRR MUR project VITALITY (ECS00000041), Spoke 2 ASTRA - Advanced Space Technologies and Research Alliance.
The authors acknowledge the support of the MUR (Italy) Department of Excellence 2023--2027 for GSSI.
Corresponding author: P. Krupa.
}%

\author[1]{Pablo~Krupa}
\ead{pkrupa@us.es}

\author[2]{Omar~Inverso}
\ead{omar.inverso@gssi.it}

\author[3]{Mirco~Tribastone}
\ead{mirco.tribastone@imtlucca.it}

\author[3]{Alberto~Bemporad}
\ead{alberto.bemporad@imtlucca.it}

\address[1]{Systems Engineering and Automation Department, Universidad de Sevilla, Seville, Spain.}
\address[2]{Gran Sasso Science Institute (GSSI), L'Aquila, Italy}
\address[3]{IMT School for Advanced Studies, Lucca, Italy}

\begin{abstract}
In safety-critical applications that rely on the solution of an optimization problem, the certification of the optimization algorithm is of vital importance. 
Certification and suboptimality results are available for a wide range of optimization algorithms.
However, a typical underlying assumption is that the operations performed by the algorithm are exact, i.e., that there is no numerical error during the mathematical operations, which is hardly a valid assumption in a real hardware implementation.
This is particularly true in the case of fixed-point hardware, where computational inaccuracies are not uncommon.
This article presents a certification procedure for the proximal gradient method for box-constrained QP problems implemented in fixed-point arithmetic.
The procedure provides a method to select the minimal fractional precision required to obtain a certain suboptimality bound, indicating the maximum number of iterations of the optimization method required to obtain it.
The procedure makes use of formal verification methods to provide arbitrarily tight bounds on the suboptimality guarantee.
We apply the proposed certification procedure on the implementation of a non-trivial model predictive controller on $32$-bit fixed-point hardware.
\end{abstract}


\end{frontmatter}
\pagestyle{fancy}
\thispagestyle{fancy}

\section{Introduction} \label{sec:intro}

Quadratic programming (QP) problems arise in various areas of systems engineering and control, such as model predictive control (MPC), see \cite{rawlings_model_2017}, or reference governors, see \cite{garone_reference_2017}, to name a few.
Various practical control-related applications, such as the ones listed above, require solving parameter-dependent QP problems at regular intervals on embedded hardware, which poses a challenge due to computational and memory limitations.
In recent years there has been a significant advance in this area due to the proposal of efficient QP solvers, some of them for generic QP problems, such as the OSQP solver presented in \cite{stellato_osqp_2020}, and some tailored to specific problems, such as the solvers proposed in \cite{krupa_implementation_2021,frison_hpipm_2020}, which address MPC optimization problems.

In many practical applications of MPC, such as safety-critical systems and space applications, the certification of the maximum number of iterations required by the optimization algorithm and a guarantee of the suboptimality of its provided solution are mandatory for real deployment.
Most solvers are based on optimization algorithms with well-known convergence and suboptimality guarantees.
The issue is that these guarantees are typically derived considering ideal conditions, e.g., under the assumption that the mathematical operations performed by the algorithm are error-free; an assumption, however, that is no longer valid when the optimization algorithm is implemented on hardware.
This is particularly noticeable on fixed-point hardware, where quantization and round-off errors may lead to significant differences with respect to the ``exact'' counterpart.
The magnitude of this difference depends on the number of fractional bits, which must be selected large enough to provide the required guarantees.

In linear-time-invariant (LTI) MPC, the use of \textit{explicit} MPC \citep{bemporad_explicit_2014} instead of an iterative solver provides a direct certification of the computation time.
However, explicit MPC may require a considerable amount of memory to implement and is only applicable to LTI systems \citep{bemporad_explicit_2014}.

In \cite{patrinos_dual_2015}, the authors present a dual gradient-projection algorithm for MPC tailored to fixed-point arithmetics.
The authors present convergence guarantees and concrete guidelines for selecting the fractional precision to obtain the required suboptimality tolerance.
The analysis is done using the notion of the \textit{inexact oracle} from \cite{devolder_first-order_2014}, which presents a generic framework for analyzing first-order optimization algorithms in which the oracle provides inexact information.
This framework can be used to derive convergence rates when inexact gradient information is available by considering the maximum error when computing the gradient.
It has, however, two downsides when applied to fixed-point arithmetic.
The first is that it only considers errors in the gradient information, i.e., it considers the other operations performed by the algorithm to be exact, which may not always be the case in fixed-point precision.
Second, the convergence results are presented in terms of the average of the iterates of the algorithm (instead of with respect to the current iterate), whose value is not generally available in fixed-point arithmetic, since its computation requires dividing by the current number of iterations.

Another approach for analyzing the error propagation is to use \textit{affine arithmetic} \citep{fang_fast_2003,vakili_finite-precision_2013}.
This framework provides less conservative error bounds than simply considering the worst-case error due to how it handles error propagation in affine operations (addition, substraction and multiplication by a constant), although the error bounds are still conservative in the presence of multiplications between variables.
The bounds can be improved by taking a probabilistic approach, as proposed in \cite{fang_fast_2003}, at the expense of no longer having a guaranteed certification.
This framework was used in \cite{nadales_efficient_2022} to certify the minimum number of fractional bits required to satisfy the desired error bound when performing Lipschitz interpolation for data-driven learning-based control.
In this case a tight certification could be proposed due to the simplicity of the algorithm, which only required affine operations.
However, its application to iterative convex optimization algorithms would provide conservative results, due, precisely, to their iterative nature along with the presence of multiplications.

Finally, formal verification methods have been used, for example,
to synthesize Lyapunov functions \citep{ahmed_Lyapunov_SMT_2020} 
or to certify the Ellipsoid method \citep{cohen_Ellipdoid_2020}, including floating-point arithmetic considerations.
\cite{simic_tight_2022} proposed a non-conservative formal verification approach for analyzing error propagation in fixed-point arithmetic. 
The technique allows to check arbitrarily tight error bounds via a bit-vector encoding into integer arithmetics and then into propositional satisfiability, which allows to use mature SAT-based technology.
The authors use their procedure to calculate accurate error bounds on a first-order optimization algorithm up to a given number of iterations.
However, due to the bit-precise encoding the computational cost of the analysis becomes prohibitive after only a few iterations, even for small optimization problems; this is an issue since first-order methods may require a significant number of iterations to converge.

In this article we analyze the implementation of the proximal gradient method (PGM) \citep{parikh_neal_proximal_2013} under fixed-point arithmetic applied to strongly convex QP problems with box constraints; a choice motivated by its simplicity but practical relevance and by its linear convergence guarantees under exact arithmetic.
We provide convergence guarantees in terms of the maximum number of iterations of the algorithm as well as suboptimality guarantees of its output.
Our certificate is based on the formal verification procedure presented in \cite{simic_tight_2022}, which we use to derive arbitrarily tight bounds on the quantities that determine the convergence and suboptimality guarantees.
The proposed approach provides a procedure for selecting the minimum fractional precision required to guarantee the given suboptimality and computation-time specifications.
Notably, the important difference with respect to \cite{simic_tight_2022} is that the verification only needs to analyze a single iteration of the PGM, thus remaining tractable for a wider range of problems.
The main features of our approach~are:
\begin{enumerate}[label=\textit{(\roman*)}, wide, labelwidth=!, labelindent=0pt, noitemsep, topsep=0pt]
    \item The subptimality guarantees are provided in terms of the output of the algorithm, instead of for the averaged iterates used in \cite{patrinos_dual_2015}.
    \item The bounds that determine the suboptimality guarantees are obtained using the formal verification procedure from \cite{simic_tight_2022}, which allows us to provide a non-conservative error-bound of the gradient computation, in that the exact maximum error committed in the computation of the gradient can be approximated to an arbitrarily large precision.
    \item The certification results are formal guarantees, instead of the probabilistic ones that would be obtained using probabilistic \textit{affine arithmetic} \citep{fang_fast_2003} or Monte Carlo analysis \citep{saracco_uncertainty_2012}.
\end{enumerate}
We present the application of our approach on a non-trivial MPC problem implemented using 32-bit fixed-point arithmetic.

\vspace*{0.3em}
\noindent\textbf{Notation:}
Given two vectors $x, y \in \R^n$, $x \leq (\geq) \; y$ denotes componentwise inequalities and $\sp{x}{y}$ is their standard inner product.
The standard Euclidean norm of a vector $x \in \R^n$ is denoted by $\|x\| \doteq \sqrt{\sp{x}{x}}$.
The closed ball of radius $r \geq 0$ in $\R^n$ is defined as the set $\B{r} \doteq \{x \in \R^n : \| x \| \leq r\}$.
The indicator function of a set $\cc{C}$ is denoted by $\Ind_{\cc{C}}$, i.e., $\Ind_{\cc{C}}(x) = 0$ if $x \in \cc{C}$ and $\Ind_{\cc{C}}(x) = \infty$ if $x \not\in \cc{C}$.
The subdifferential of a function $f$ is denoted by $\partial f$.
We denote by $\Rpp$ the set of strictly positive real numbers.
We denote by $\R^n\pq \subset \R^n$ the set of vectors whose integer and fractional parts are representable using $p$ and $q$ binary digits, respectively.
This notion readily extends to the space of matrices $\R^{n \times m}$, where $\R^{n \times m}\pq \subset \R^{n \times m}$ represents the space of matrices whose every element is representable in $\R\pq$.
For any set $\cc{C} \subseteq \R^n$, we denote $\cc{C}\pq = \set{x \in \R^n\pq}{x \in \cc{C}}$.
It is obvious that $x \in \cc{C}\pq \implies x \in \cc{C}$, but not vice-versa.
It is also easy to see that $\cc{C}\pq[p'.q'] \subseteq \cc{C}\pq$ for any $p \geq p'$ and $q \geq q'$, and therefore that $x \in \cc{C}\pq[p'.q'] \implies x \in \cc{C}\pq$.
\vspace*{-0.3em}

\section{Exact proximal gradient method} \label{sec:PGMexact}

We consider the class of strongly-convex QP problems
\begin{equation*}
    \begin{aligned}
        \min\limits_{x \in \R^n} &\; \frac{1}{2}x\T Q x + c\T x\\
        \st &\; \ell \leq x \leq u,
    \end{aligned}
    \label{eq:QPbox}
    \tag{$\pP$}
\end{equation*}
where $c, \ell, u, \in \R^n_{(p'.q')}$, with $\ell \leq u$, and $Q \in \R^{n \times n}_{(p'.q')}$, for some finite positive $p'$ and $q'$.
In this article we are interested in finding a suboptimal solution of problem~\eqref{eq:QPbox}, using the proximal gradient method (PGM) \citep{parikh_neal_proximal_2013}, for any given realization of the ingredients of problem \eqref{eq:QPbox} satisfying
\begin{subequations} \label{eq:ingredients}
\begin{align}
    Q &\in \cQ \subset \R^{n \times n}_{(p'.q')}, \; c \in \set{c \in \R^n_{(p'.q')}}{c_{\tmin} \leq c \leq c_{\tmax}}, \\
    \ell &\in \set{\ell \in \R^n_{(p'.q')}}{\ell_{\tmin} \leq \ell \leq \ell_{\tmax}}, \\
    u &\in \set{u \in \R^n_{(p'.q')}}{u_{\tmin} \leq u \leq u_{\tmax}},
\end{align}
\end{subequations}
where $c_{\tmin}, c_{\tmax}, \ell_{\tmin}, \ell_{\tmax}, u_{\tmin}, u_{\tmax} \in \R^n_{(p'.q')}$, $c_{\tmin} \leq c_{\tmax}$, $\ell_{\tmin} \leq \ell_{\tmax}$, $u_{\tmin} \leq u_{\tmax}$, and $\cQ$ is a compact set.
This problem arises in numerous control-related fields, such as MPC problems that can be written in compact form \citep{richter_computational_2012}; it can also be applied to a dual setting provided that a compact bound on the dual variables is available \citep{patrinos_dual_2015}.
In particular, we are interested in providing convergence guarantees when implementing the PGM in fixed-point arithmetic.
To that end, let us start by recalling the classical ``exact'' PGM, i.e., its implementation when operating using exact arithmetic.

We denote by $\cP$ the set of problems \eqref{eq:QPbox} whose ingredients satisfy \eqref{eq:ingredients} and $\ell \leq u$.
Let $\cc{X} \doteq \set{x \in \R^n}{\ell \leq x \leq u}$, $f: \R^n \rightarrow \R$ be given by $f(x) = \frac{1}{2} x\T Q x + c\T x$, $x^* \in \R^n$ be the optimal solution of problem \eqref{eq:QPbox} and $f^*$ its optimal value, i.e., ${f^* = f(x^*)}$.
Let $\L, \ss \in \Rpp$ be the largest and smallest smoothness and strong convexity parameters of $f$ for any realization of \eqref{eq:QPbox}$\,\in \cP$, i.e., the scalars for which the well-known inequalities
\vspace*{-0.5em}
\begin{align*}
    f(x) &\leq f(y) + \sp{ \g{y} }{ x - y} + \frac{\L}{2} \| x - y \|^2, \\
    f(x) &\geq f(y) + \sp{ \g{y} }{ x - y} + \frac{\ss}{2} \| x - y \|^2,
\end{align*}
are satisfied $\forall x, y \in \R^n$ for any given \eqref{eq:QPbox}$\,\in \cP$.

\begin{algorithm}[t]
    \DontPrintSemicolon
    \caption{PGM applied to \eqref{eq:QPbox}$\,\in\cP$} \label{alg:PGMexact}
    \KwIn{$x^0 \in \cX^0$, $0 < \t \leq \L^{-1}$}
    \MyRepeatEver{}{
        $x^{k+1} \gets \CGM[\t](x^k) = \tmin\{u, \tmax\{\ell, x^k - \t \nabla f (x^k) \}\}$\; \label{step:PGMexact}
    }
\end{algorithm}

Algorithm~\ref{alg:PGMexact} shows the PGM applied to a realization of problem \eqref{eq:QPbox}$\,\in\cP$, where the $\tmin$ and $\tmax$ operators are taken componentwise and $0 < \t \leq \L^{-1}$.
It generates a sequence $\{x^k\} \in \cX$ starting at an initial point $x^0 \in \Xinit$, where $\Xinit \subseteq \cX$ is the set of possible initial guesses.
Step~\ref{step:PGMexact} of the algorithm performs the operator $\CGM: \cX \rightarrow \cX$ given by
\begin{equation*}
    \CGM[\t](x) \doteq \arg\min\limits_{y \in \cX}  \frac{1}{2} \left\| y - \left( x - \t \g{x} \right) \right\|^2,
\end{equation*}
which, when particularized to problem \eqref{eq:QPbox}, is the evaluation of the \textit{proximal} operator \citep{parikh_neal_proximal_2013} at $x - \t \g{x}$.
Operator $\CGM[\t]$ in this setting can also be viewed as the evaluation of the composite gradient mapping \citep[\S 2]{nesterov_gradient_2013}.
Furthermore, Algorithm~\ref{alg:PGMexact} is equivalent in this case to the projected gradient method \citep[\S 4,2]{parikh_neal_proximal_2013}.

The following theorem recalls the linear convergence of Algorithm~\ref{alg:PGMexact} when working under ``exact'' arithmetic.
In the following section we will derive a similar result when operating under fixed-point arithmetic.

\begin{theorem}[Theorem 10.29 from \cite{beck_first-order_2017}] \label{theo:convergence:exact}
Let $\{x^k\}$ be the sequence generated by Algorithm~\ref{alg:PGMexact} starting at $x^0 \in \Xinit$ applied to \eqref{eq:QPbox}$\,\in\cP$.
Then,
\begin{Rlist}
    \item $ \| x^k - x^* \|^2 \leq \left( 1 - \t \ss \right)^{k} \| x^0 - x^* \|^2, \; \forall k \geq 0$,
    \item $ f(x^{k}) - f^* \leq \fracg{1}{2 \t} \left( 1 - \t \ss \right)^{k} \| x^0 - x^* \|^2, \; \forall k \geq 1$.
\end{Rlist}
\end{theorem}

\section{Proximal gradient method in fixed-point arithmetic} \label{sec:PGMfp}

The maximum number of iterations of Algorithm~\ref{alg:PGMexact} required to guarantee a given suboptimality can be certified using Theorem~\ref{theo:convergence:exact}.
Indeed, an immediate result of Theorem~\ref{theo:convergence:exact} is that $\| x^k - x^* \|^2 \leq \epsilon$, for a given $\epsilon \in \Rpp$, is satisfied for every iteration $k$ satisfying
\begin{equation} \label{eq:k:exact}
    k \geq \frac{\log \left( \fracg{\epsilon}{\| x^0 - x^* \|^2} \right)}{\log(1 - \t \ss)}.
\end{equation}

We are now interested in providing an iteration and suboptimality certification when Algorithm~\ref{alg:PGMexact} is implemented using fixed-point arithmetic.
Therefore, let us consider Algorithm~\ref{alg:PGMexact} when working under fixed-point arithmetic for some predetermined choice of integer and fractional precision $(p.q)$ satisfying $p \geq p'$ and $q \geq q'$, where we recall that $(p'.q')$ is the precision under which the ingredients of \eqref{eq:QPbox} are representable.
Under this paradigm, we can view the implementation of the PGM algorithm as performing an inexact proximal operator, where the source of inexactness is due to the fixed-point arithmetic and representation of variables.

\begin{algorithm}[t]
    \DontPrintSemicolon
    \caption{PGM applied to \eqref{eq:QPbox}$\,\in\cP$ under fixed-point arithmetic with precision $(p.q)$} \label{alg:PGMfp}
    \KwIn{$\fpx^0 \in \fpXinit\pq$, $0 < \t \leq \L^{-1}$, $\fp{\epsilon} \in \Rpp\pq$, $k_\tmax \in \Rpp$}
    \MyRepeat{$\distfp < \fp{\epsilon}$ or $k \geq k_\tmax$ \label{step:PGMfp_exit}}{
        $\fpx^{k+1} \gets \min\{u, \max\{\ell, \fpx^k - \gfp[\fpx^k] \}\}$\; \label{step:PGMfp}
    }
\end{algorithm}

Algorithm~\ref{alg:PGMfp} shows the implementation of Algorithm~\ref{alg:PGMexact} under fixed-point arithmetic.
It generates a sequence of iterates $\{\fpx^k\} \in \cX\pq$ starting from an initial point $\fpx^0 \in \fpXinit\pq \subseteq \cX\pq$.
Step~\ref{step:PGMfp}  evaluates the operator $\fp{g}_\t : \cX\pq \rightarrow \cX\pq$, which is defined as the operator that performs the computation of $\t \g{\cdot}$ in the fixed-point paradigm.
That is, $\gfp[\fpx^k]$ returns the result of the evaluation of $\t \g{\fpx^k}$ when performed under fixed-point arithmetic.
Thus, $\t \in \Rpp\pq$ is an obvious requirement of Algorithm~\ref{alg:PGMfp}.
Additionally, the algorithm includes an exit condition given by the satisfaction of the condition $\distfp < \fp{\epsilon}$, where $\fp{d}^2 : \cX\pq \rightarrow \Rpp\pq \cup \{0\}$ is the operator that performs the computation of $\| \fpx^{k+1} - \fpx^{k} \|^2$ in fixed-point arithmetic.
This exit condition plays a key role in the suboptimality guarantees provided in this section.

The value of $\gfp[\fpx^k]$ will generally differ from the exact value of the expression $\t \g{\fpx^k}$ due to the arithmetic errors that occur when using fixed-point arithmetic, thus leading to the source of discrepancy between the sequences generated by Algorithms~\ref{alg:PGMexact} and \ref{alg:PGMfp}.
The magnitude of this discrepancy, which we formalize in the following definition, will depend on the value of the fractional precision $q$, with higher values of $q$ obviously leading to smaller errors.

\begin{definition} \label{def:Omega}
    Given a choice of $q$, we denote by $\Omega \in \Rpp$ a scalar satisfying $\| \t \g{\fpx} - \gfp[\fpx] \| \leq \Omega$, $\forall \fpx \in \cX\pq$, $\forall$\eqref{eq:QPbox}$\,\in \cP$.
\end{definition}

The rest of the operations in Step~\ref{step:PGMfp} of Algorithm~\ref{alg:PGMfp} do not incur any additional error, as formally stated in the following lemma, so they do not contribute towards the discrepancy between the ``exact'' and ``fixed-point'' implementations.

\begin{lemma} \label{lemma:exact:projection}
    Let $\underline{v}, \overline{v}, \hat{x}, \hat{y} \in \R^n\pq$, with $\underline{v} \leq \overline{v}$, and consider the set $\cc{C} \doteq \set{v \in \R^n}{\underline{v} \leq v \leq \overline{v}}$.
    Then, the result of the fixed-point computations $\min\{ \overline{v},\, \max\{ \underline{v},\, \hat{x} + \hat{y} \} \}$ performed with any precision $(\hat{p}.\hat{q})$ satisfying $\hat{p} \geq p$ and $\hat{q} \geq q$ is the exact Euclidean projection of $\hat{x} + \hat{y}$ onto $\cc{C}$ if $\hat{x} + \hat{y}$ does not result in an overflow.
\end{lemma}

The proof of the lemma is omitted because the claim is a direct result of the fact that the $\min$, $\max$ and addition operations do not incur in any error under fixed-point arithmetic as long as there is no overflow in the addition.

\begin{corollary} \label{corr:exact:projection}
    Variable $\fpx^{k+1}$ obtained from Step~\ref{step:PGMfp} of Algorithm~\ref{alg:PGMfp} is the exact Euclidean projection of $\fpx^k - \gfp$ onto $\cX$, assuming that no overflow occurs during the computations.
    Therefore, $\fpx^k \in \cX\pq \subset \cX$, $\forall k \geq 0$.
\end{corollary}

The reader will note that Lemma~\ref{lemma:exact:projection} is only applicable to Algorithm~\ref{alg:PGMfp} as long as there is no overflow during its execution, i.e., if the integer precision $p$ is large enough.
The certification tool presented in Section~\ref{sec:verification} can be used to compute the minimum value of $p$ required to avoid overflow.
Thus, we henceforth simply consider that $p$ is chosen so that no overflow occurs during the execution of Algorithm~\ref{alg:PGMfp}.

Another useful consequence of Lemma~\ref{lemma:exact:projection} is presented in the following lemma, which states that a scalar $\Omega$ satisfying Definition~\ref{def:Omega} also bounds the error in the computation of $\fpx^{k+1}$.

\begin{lemma} \label{lemma:bound:err:x}
Consider Algorithm~\ref{alg:PGMfp} and let $\Omega$ satisfy Definition~\ref{def:Omega}.
Then, $\| \fpx^{k+1} - \CGM(\fpx^k) \| \leq \Omega$, $\forall \fpx^k \in \cX\pq$, $\forall$\eqref{eq:QPbox}$\,\in \cP$.
\end{lemma}

\begin{proof}
The claim is a direct consequence of Corollary~\ref{corr:exact:projection} and the fact that the projection operator to non-empty closed convex sets is non-expansive~\cite[\S 3.1]{ryu_primer_2016}.
\end{proof}

We now present the main result of this section, where we characterize the local linear convergence of Algorithm~\ref{alg:PGMfp}, in terms of the error-bound $\Omega$, when sufficiently far away from the optimal solution.

\begin{theorem} \label{theo:convergence:fp}
Let $\{\fpx^k\}$ be the sequence generated by Algorithm~\ref{alg:PGMfp} applied to a realization of problem \eqref{eq:QPbox}$\,\in\cP$ with starting point $\fpx^0 \in \fpXinit\pq$ and taking $0 < \t \leq \L^{-1}$, $\t \in \R\pq$.
Choose $\epsilon \in \R^+$ satisfying $\epsilon \t \ss > 4 \Omega$, where $\Omega \in \R^+$ is given by Definition~\ref{def:Omega}.
Then, as long as $\| \fpx^{k+1} - x^* \| \geq \fracg{\epsilon}{2}$, the sequence $\{\fpx^k\}$ satisfies:
    \begin{RlistL}
    \item $\| \fpx^k - x^* \|^2 \leq \left( \fracg{1 - \t \ss }{1 - 4 \Omega \epsilon^{-1}} \right)^{k} \| \fpx^0 - x^* \|^2, \; \forall k \geq 0$. \label{theo:convergence:fp:xk}
    \item $f(\fpx^{k}) {-} f^* \leq \fracg{1 {-} 4 \Omega \epsilon^{-1}}{2 \t} \left( \fracg{1 - \t \ss }{1 {-} 4 \Omega \epsilon^{-1}} \right)^{k} \| \fpx^0 - x^* \|^2, \; \forall k \geq 1$. \label{theo:convergence:fp:f}
    \end{RlistL}
\end{theorem}

\begin{proof}\renewcommand{\qedsymbol}{}
See Appendix \ref{app:proofs}.
\end{proof}

Theorem~\ref{theo:convergence:fp} provides a linear convergence result similar to the one shown in Theorem~\ref{theo:convergence:exact}, but where the convergence constant degrades by the factor $(1 - 4 \Omega \epsilon^{-1})$.
That is, the convergence guarantee worsens as $\Omega$ increases and as the desired suboptimality tolerance $\epsilon$ decreases.
Since the theorem only holds as long as $\| \fpx^{k+1} - x^* \| \geq \fracg{\epsilon}{2}$, we need to be able to check the satisfaction of this condition during the execution of Algorithm~\ref{alg:PGMfp}.

\begin{lemma} \label{lemma:dist_opt}
    $\| \fpx - \CGM[\t](\fpx) \| \leq 2 \| \fpx - x^* \|$, $\forall \fpx \in \cX$.
\end{lemma}

\begin{proof}
From \cite[Property 1.(i)]{alamo_restart_2019}, particularized to our problem formulation and notation, we have that
\begin{equation*}
    f(\CGM[\t](\fpx)) - f^* \leq \ti \sp{\fpx - \CGM[\t](\fpx)}{\fpx - x^*} - \frac{1}{2\t} \| \fpx - \CGM[\t](\fpx) \|^2,
\end{equation*}
which along with $f(\CGM[\rho](\fpx)) - f^* \geq 0$, leads to
\begin{equation*}
    \frac{1}{2} \| \fpx - \CGM[\t](\fpx) \|^2 \leq \sp{\fpx - \CGM[\t](\fpx)}{\fpx - x^*} \leq \| \fpx - \CGM[\t](\fpx) \| \cdot \| \fpx - x^* \|
\end{equation*}
by making use of the Cauchy-Schwarz inequality.\qedhere
\end{proof}

The previous lemma allows us to guarantee that the condition $\| \fpx^{k} - x^* \| \geq \fracg{\epsilon}{2}$ in Theorem~\ref{theo:convergence:fp} holds for the iterates of Algorithm~\ref{alg:PGMfp} as long as we can guarantee that $\| \fpx^k - \CGM[\t](\fpx^k) \| \geq \epsilon$.
The following assumption allows us to use the exit condition of Algorithm~\ref{alg:PGMfp} as a means to guarantee that the condition $\| \fpx^k - x^* \| \geq \fracg{\epsilon}{2}$ is satisfied at iteration $k$.

\begin{assumption} \label{ass:epsilon_hat}
The exit tolerance $\hat{\epsilon}$ of Algorithm~\ref{alg:PGMfp} satisfies $\distfp \geq \hat{\epsilon} \implies \| \fpx^k - \CGM[\rho](\fpx^k) \|^2 \geq \epsilon^2$.
\end{assumption}

In the following section we present a tractable procedure for certifying the satisfaction of this assumption.
In practice, we find that one can choose $\epsilon$ so that the smallest value of $\fp{\epsilon}$ satisfying Assumption~\ref{ass:epsilon_hat} is the smallest positive representable number in precision $q$, i.e., $2^{-q}$.
In this case the exit condition of Algorithm~\ref{alg:PGMfp} becomes $\distfp = 0$.

Under Assumption~\ref{ass:epsilon_hat}, the convergence guarantee provided in Theorem~\ref{theo:convergence:fp} holds as long as the exit condition $\distfp < \fp{\epsilon}$ is not satisfied.
If it is satisfied at some iteration $k$, then the convergence guarantee provided by Theorem~\ref{theo:convergence:fp} is only guaranteed to hold until iteration $k-1$.
The following theorem provides the suboptimality guarantees of Algorithm~\ref{alg:PGMfp} when the exit condition is satisfied at some iteration $k$.
The result makes use of the bounds provided in the following definition.
The following section will provide computationally tractable procedures for computing arbitrarily tight values of said bounds.

\begin{definition} \label{def:exit:bounds}
Consider Algorithm~\ref{alg:PGMfp}.
For a given choice of fractional precision $q$ and tolerance $\fp{\epsilon}$, we denote by $\delta, \omega, \Theta \in \Rpp$ the scalars satisfying
\begin{equation*}
    \| \fpx^k - \CGM[\t](\fpx^k) \| \leq \delta, \; \| \fpx^{k+1} -  \CGM[\t](\fpx^k) \| \leq \omega, \; \| \fpx^{k} -  \fpx^{k+1} \| \leq \Theta
\end{equation*}
for all $\fpx^k \in \cX\pq$ satisfying $\distfp < \fp{\epsilon}$, $\forall$\eqref{eq:QPbox}$\,\in\cP$.
\end{definition}

\begin{theorem} \label{theo:convergence:exit}
Let $\{\fpx^k\}$ be the sequence generated by Algorithm~\ref{alg:PGMfp} applied to a realization of problem \eqref{eq:QPbox}$\,\in\cP$ with starting point $\fpx^0 \in \fpXinit\pq$ and taking $0 < \t \leq \L^{-1}$, $\t \in \R\pq$.
Denote $T \doteq \ss^{-1} (\ti + \L)$.
Then, if $\distfp[\fpx^k] < \fp{\epsilon}$,
\begin{Rlist}
\item $\| \fpx^{k+1} - x^* \| \leq \omega + \delta T$, \label{theo:convergence:exit:x}
    \item $f(\fpx^{k+1}) - f^* \leq \ti \left( (\Theta + \Omega) (\omega + \delta T) + \fracg{1}{2} \Theta^2 \right)$. \label{theo:convergence:exit:f}
\end{Rlist}
\end{theorem}

\begin{proof}
    From \cite[Lemma 3]{nesterov_gradient_2013}, particularized to our notation, we have that $\| \fpx^k - \CGM[\t](\fpx^k) \| \geq T^{-1} \| \CGM[\t](\fpx^k) - x^* \|$, which leads to
    $\| \fpx^k - \CGM[\t](\fpx^k) \| \leq \delta {\implies} \| \CGM[\t](\fpx^k) - x^* \| \leq T \delta$.
    Claim~\ref{theo:convergence:exit:x} follows from adding the previous inequality with $\| \fpx^{k+1} - \CGM[\t](\fpx^k) \| \leq \omega$ and applying the triangle inequality.
    By the same procedure, we also derive $\| \fpx^k - x^* \| \leq \delta ( T + 1)$.
    Claim~\ref{theo:convergence:exit:f} then follows from particularizing Lemma~\ref{lemma:PGMfp:distf}.\eqref{eq:PGMfp:distf:1} to $y = x^*$, using the Cauchy-Schwarz inequality and then taking the previous inequalities along with the inequalities presented in Definitions~\ref{def:Omega} and \ref{def:exit:bounds}.
\end{proof}

The following corollary gathers the guarantees that are obtained from Algorithm~\ref{alg:PGMfp} in terms of the error-bounds and tolerances presented throughout this section.

\begin{corollary}[Suboptimality guarantee of Algorithm~\ref{alg:PGMfp}] \label{corr:certification}
Let $k_{\rm max} \geq \log\left( \frac{\epsilon^2}{4 D} \right)/\log(C)$, where
$C \doteq \frac{1 - \t \ss}{1 - 4 \Omega \epsilon^{-1}}$ and $D$  satisfies
$D \geq \max_{\fpx \in \fpXinit\pq} \| \fpx - x^* \|^2$.
The following hold:
\begin{Rlist}
\item If $\distfp \geq \fp{\epsilon}$ for all $k = \{0, \dots k_{\tmax} \}$, then
    \begin{equation*}
        \| \fpx^{k_{\tmax}} - x^* \| \leq \fracg{\epsilon}{2} \;\; \text{and} \;\; f(\fpx^{k_{\tmax}}) - f^* \leq \frac{(\epsilon^2 - 4 \Omega \epsilon)}{8 \t}.
    \end{equation*} \label{corr:certification:kmax}
\item If $\distfp < \fp{\epsilon}$ then $\| \fpx^{k+1} - x^* \| \leq \omega + \delta T$ and
    \begin{equation*}
        f(\fpx^{k+1}) - f^* \leq \ti \left( (\Theta + \Omega) (\omega + \delta T) + \frac{1}{2} \Theta^2 \right).
    \end{equation*} \label{corr:certification:dist}
\end{Rlist}
\end{corollary}

\begin{remark} \label{rem:compute:omega}
From Lemma~\ref{lemma:bound:err:x} we have that the bound $\omega$ from Definition~\ref{def:exit:bounds} can be substituted by $\Omega$ from Definition~\ref{def:Omega}.
We find that for high fractional precision there can be an insignificant difference between the two quantities, thus not meriting the additional computation time required to compute $\omega$.
\end{remark}

\section{Obtaining error-bounds for Algorithm \ref{alg:PGMexact}} \label{sec:verification}

This section presents procedures for obtaining arbitrarily tight values of the bounds $\Omega$, $\delta$, $\omega$ and $\Theta$ introduced in Definitions~\ref{def:Omega} and \ref{def:exit:bounds} as well as a procedure for checking the satisfaction of Assumption~\ref{ass:epsilon_hat}.
The procedures consist on the application of the formal verification procedure presented in \cite{simic_tight_2022} to check the conditions provided in the definitions and assumption.
We first introduce the formal analysis technique for fixed-point arithmetic, and then show how it fits within our certification process.

\subsection{Formal verification for fixed-point arithmetic} \label{sec:verification:cseq}

\cite{simic_tight_2022} consider the problem of estimating the numerical accuracy of algorithms in fixed-point arithmetic with variables of arbitrary precision and possibly non-deterministic values.
The idea is to re-compute in a greater precision the result of each fixed-point operation, so that the numerical error can be estimated based on the difference between the two values; at the same time, the different errors are in turn accounted for and propagated through the re-computations.
When sufficient precision is used to store the re-computed values, this yields an accurate error tracking for each variable at any point of the computation.
The techniques relies on a bit-precise encoding to transform the sequences of operations under analysis into operations in integer arithmetic over vectors of bits; these are in turn encoded as a SAT formula that is satisfiable if and only if the algorithm under analysis exceeds a given bound on the numerical error. 
The technique is quite accurate in that it allows to formally verify arbitrarily tight bounds on the numerical error up to a given number of iterations.
On the other hand, the required program unfolding pass along with the bit-vector encoding ends up introducing considerable overhead; the analysis can become quickly intractable, even for small problems and a few iterations.
We note that this is not a problem in our case, as the algorithms that we need to analyze have been formulated so as to only require a single iteration.

\begin{algorithm}[t]
    \DontPrintSemicolon
    \caption{Example of input algorithm for formal verification} \label{alg:cseq:example}
    \SetAlCapSkip{1em}
    \KwParameters{$(p.q)$, $b \in \R\pq$, $\hat{a} \in \Rpp$, $\xi \in \R$, $\chi \in \Rpp$}
    \KwNonDeterministic{$a \in \mathcal{A} = \set{a \in \R^m}{|| a ||_\infty \leq \hat{a}}$}
    $r = \sp{a}{a}$\;
    $\mu = b r$\;
    \lIf{$ | \textnormal{\texttt{err}}(r) | > \chi $}{\texttt{FAIL}} \label{alg:cseq:example:err}
    \lIf{$ \textnormal{\texttt{exact}}(\mu) < \xi $}{\texttt{FAIL}} \label{alg:cseq:example:exact}
\end{algorithm}

We now provide an illustrative example to introduce the concepts and notation relevant to this article.
Consider Algorithm~\ref{alg:cseq:example} working under a given fixed-point precision $(p.q)$, where $r$ is a fixed-point variable, $\exact{r}$ is its ``exact'' counterpart and $\err{r}$ is its error, i.e., $\exact{r} = r + \err{r}$.
The same applies for the other variables, such as $\mu$, whose exact value will generally differ from its value computed in fixed-point arithmetic due to the multiplication operations.
The error of $r$ is propagated when computing the error of $\mu$, i.e., $\exact{\mu}$ will contain the value of $\sp{a}{a}$.
We can perform assessments on the fixed-point variables, their errors and their exact values (as shown in Steps~\ref{alg:cseq:example:err} and \ref{alg:cseq:example:exact}) for all possible values of the non-deterministic inputs (variable $a$ in this case).
If all the assessments are satisfied for all possible values of the non-deterministic variables then the procedure will return a \texttt{PASS}.
Otherwise, it will return a \texttt{FAIL}.

\begin{algorithm}[t!]
    \DontPrintSemicolon
    \caption{Algorithm for asserting Assumption~\ref{ass:epsilon_hat}} \label{alg:verify:ass:eps}
    \SetAlCapSkip{1em}
    \KwParameters{$(p.q)$, $\fp{\epsilon}$, $\epsilon^2$}
    \KwNonDeterministic{$\fpx^0 \in \cX\pq$, \eqref{eq:QPbox}$\,\in\cP$}
    $\fpx^{k+1} \gets \min\{u, \max\{\ell, \fpx^k - \gfp[\fpx^k] \}\}$\;
    $d \gets \sp{\fpx^{k} - \fpx^{k+1}}{\fpx^{k} - \fpx^{k+1}}$\;
    \lIf{$d \geq \hat{\epsilon}$ \KwAnd $\exact{d} < \epsilon^2$}{\texttt{FAIL}}
\end{algorithm}

\begin{algorithm}[t!]
    \DontPrintSemicolon
    \caption{Algorithm for deriving $\Omega$} \label{alg:verify:Omega}
    \SetAlCapSkip{1em}
    \KwParameters{$(p.q)$, $\Omega^2$}
    \KwNonDeterministic{$\fpx^0 \in \cX\pq$, \eqref{eq:QPbox}$\,\in\cP$}
    $\hat{g} \gets \t (Q \fpx^0 + c)$\;
    $e \gets \err{\hat{g}}$\;
    $v \gets \sp{e}{e}$\;
    \lIf{$ \exact{v} > \Omega^2$}{\texttt{FAIL}} \label{step:verify:Omega:FAIL}
\end{algorithm}

\begin{algorithm}[t!]
    \DontPrintSemicolon
    \caption{Algorithm for deriving Def.~\ref{def:exit:bounds} bounds} \label{alg:verify:exit}
    \SetAlCapSkip{1em}
    \KwSelect{b $\in \{ \delta^2, \omega^2, \Theta^2 \}$}
    \KwParameters{$(p.q)$, $\fp{\epsilon}$}
    \KwNonDeterministic{$\fpx^0 \in \cX\pq$, \eqref{eq:QPbox}$\,\in\cP$}
    $\fpx^{k+1} \gets \min\{u, \max\{\ell, \fpx^k - \gfp[\fpx^k] \}\}$\;
    \If{$\sp{\fpx^{k} - \fpx^{k+1}}{\fpx^{k} - \fpx^{k+1}} < \hat{\epsilon}$}{
        \Switch{\normalfont{b}}{
            \lCase{$\delta^2$}{$s \gets \fpx^k - \fpx^{k+1}$}
            \lCase{$\omega^2$}{$s \gets \err{\fpx^{k+1}}$}
            \Case{$\Theta^2$}{
                $\err{\fpx^{k+1}} \gets 0$\;
                $s \gets \fpx^{k} - \fpx^{k+1}$\;
            }
        }
    }
    $v \gets \sp{s}{s}$\;
    \lIf{$ \exact{v} > \normalfont{b}$}{\texttt{FAIL}}
\end{algorithm}

\begin{example} \label{example:cseq}
We run the verification procedure presented in \cite{simic_tight_2022} on Algorithm~\ref{alg:cseq:example} with $p = 8$, $q = 8$, $m = 20$, $\hat{a} = 0.125$, $b = 1.5$, $\xi = 0$, and $\chi = 0.069580078125$.
We obtain a {\normalfont\texttt{PASS}}, that is, there is no value of $a$ inside the box defined by $\hat{a}$ for which the assertions stated in Steps~\ref{alg:cseq:example:err} and \ref{alg:cseq:example:exact} of Algorithm~\ref{alg:cseq:example} are violated.
The result of this example highlights one of the main benefit of this procedure, which is the bound $\chi = 0.069580078125$.
Variable $r$ is the inner product of $a$ by itself.
As stated in \cite{patrinos_dual_2015}, the standard theoretical bound for the maximum error committed by an inner product is given by $\texttt{err}(\sp{a}{a}) \leq 2^{-q} m$, which is equal to $0.078125$ for the values of $q$ and $m$ in this example.
However, the formal verification procedure has found that this theoretical bound can be improved to $\texttt{err}(\sp{a}{a}) = 0.069580078125$, which is a $12.28\%$ improvement.
That is, the procedure may lead to tighter bounds on the errors committed by the fixed-point algorithm than the ones obtained by simply evaluating its execution under the theoretical worst-case scenario.
The procedure finishes after $\sim5$s using a standard machine with an Intel i5 processor.
\end{example}

\subsection{Verification procedure for Algorithm~\ref{alg:PGMfp}} \label{sec:verification:bounds}

The application of the verification tool presented in the previous subsection to Algorithm~\ref{alg:verify:ass:eps} certifies if the given $\fp{\epsilon}$ and $\epsilon^2$ satisfy Assumption~\ref{ass:epsilon_hat}.
Its application to Algorithms~\ref{alg:verify:Omega} certifies if the provided value of $\Omega$ satisfies the condition presented in Definition~\ref{def:Omega}.
A \texttt{PASS} indicates that the given $\Omega^2$ satisfies the condition $\| \CGM[\t](\fpx^k) - \fpx^{k+1} \|^2 \leq \Omega^2$, $\forall \fpx^k \in \cX\pq$, $\forall$\eqref{eq:QPbox}$\,\in\cP$, whereas a \texttt{FAIL} indicates that there is at least one combination of $\fpx^k \in \cX\pq$ and \eqref{eq:QPbox}$\,\in\cP$ for which the condition is not satisfied.
Note that the condition is asserted with respect to $\Omega^2$, since the verification tool does not allow the use of the square-root operation.
Similarly, for a given value of $\fp{\epsilon}$, Algorithm~\ref{alg:verify:exit} is used to certify if the bounds $\delta$, $\omega$ or $\Theta$ satisfy the conditions presented in Definition~\ref{def:exit:bounds}.
Arbitrarily tight values of $\Omega$, $\delta$, $\omega$ or $\Theta$ can be obtained by applying the bisection method to Algorithms~\ref{alg:verify:Omega} and \ref{alg:verify:exit}.
Note that Algorithms~\ref{alg:verify:ass:eps}, \ref{alg:verify:Omega} and \ref{alg:verify:exit} only execute a single iteration of Algorithm~\ref{alg:PGMfp}.
Thus, the proposed verification procedure remains tractable for moderately-sized problems, in contrast with the approach taken in \cite{simic_tight_2022}.

\begin{remark}
The verification tool from \cite{simic_tight_2022} can be configured to return a \texttt{FAIL} in the event of a numerical overflow.
Thus, we can verify that no overflow occurs in Algorithm~\ref{alg:PGMfp} for a given choice of the integer precision $p$ if the tool does not fail due to an overflow when applied to Algorithm~\ref{alg:verify:ass:eps}.
\end{remark}

\vspace*{-1.0em}
\section{Numerical case study} \label{sec:results}

We apply the verification procedures presented in the previous section to certify the fixed-point implementation of the PGM to solve the optimization problem of a linear MPC controller for a discrete-time, time-invariant system given by a state-space model $\tx(t) = A \tx(t) + B \tu(t)$, where $\tx(t) \in \R^{n_x}$ and $\tu(t) \in \R^{n_u}$ are the state and control input at sample time~$t$.

In particular, we consider the system of three masses connected by springs presented in \cite[\S 3]{krupa_sparse_2021}, where we take the mass all three objects equal to $1$kg and the spring constants as $1$N/m.
The $6$-dimensional system state is given by the position and velocity of each of the three objects, while the control input is given by the two external forces applied to the outer objects.
We take the following MPC formulation:
\begin{subequations} \label{eq:MPC} 
\begin{align}  
    \min \;& \Sum{i = 0}{N_p-1} \left( \| \tx_i {-} \tx_r \|^2_{W_x} {+} \| \tu_i {-} \tu_r \|^2_{W_u} \right) + \| \tx_{N_p} - \tx_r \|^2_P \\
    \st& \; \tx_0 = \tx(t) \label{eq:MPC:initial} \\
        & \; \tx_{i+1} = A \tx_i + B \tu_i, \; \forall i\in\Z_0^{N_p - 1} \label{eq:MPC:prediction} \\
        & \; \tu_{i} = \tu_{N_c-1}, \; \forall i \in \Z_{N_c}^{N_p - 1} \\
        & \; \tu_{-} \leq \tu_i \leq \tu_{+}, \; \forall i\in\Z_0^{N_c-1},
\end{align}
\end{subequations}
where $N_c \in \Rpp$ is the control horizon; $N_p \geq N_c$ is the prediction horizon; $W_x, P \in \R^{n_x \times n_x}$ and $W_u \in \R^{n_u \times n_u}$ are positive definite; $(x_r, u_r)$ are the state and input references; and $\tu_{-}, \tu_{+} \in \R^{n_u}$ satisfying $\tu_{-} \leq \tu_{+}$ define the bounds on the control input.
We take $N_c = 2$, $N_p = 5$, $\tu_{+} = (0.5, 0.5)$, $\tu_{-} = -\tu_{+}$, $W_x = 0.5 I_{n_x}$, $W_u = 0.25 I_{n_u}$ and $P$ as the solution of the associated discrete algebraic Riccati equation.
Problem \eqref{eq:MPC} can be transformed into \eqref{eq:QPbox} by eliminating the states and rewriting it in condensed form, see e.g., \cite{richter_computational_2012, jerez_condensed_2011}, leading to a $(n_u N_c)$-dimensional QP problem.
In this case, ingredients $Q$, $\ell$ and $u$ of problem \eqref{eq:QPbox} are fixed, whereas the value of $c$ will depend on the value of the reference $(x_r, u_r)$ as well as the current state $\tx(t)$.
We compute $c_\tmin$ and $c_\tmax$ by assuming that the position of the objects belong to the interval $[-0.5, 0.5]$m and the velocities to $[-1, 1]$m/s.
A non-deterministic $Q$ would be taken if we allowed the possibility of changing the weight $W_u$ online or if we considered a time-varying model of the system.

We now certify the PGM applied to the resulting condensed MPC problem when implemented in fixed-point arithmetic on a $32$-bit device, where we take $p = 10$ for the integer precision and $q = 21$ for the fractional precision (the remaining bit is used for storing the sign).
We store the matrices of the QP problem in the selected precision.
The resulting problem has $\L = 4.9645$ and $\ss = 0.3532$.
We take $\rho$ as the largest number representable in $\R\pq$ that satisfies $\rho \leq 1/\L$.

{\renewcommand{\arraystretch}{1.1}%
\begin{table*}[t]
    \centering
    \setlength\tabcolsep{5pt} 
        \begin{tabular}{|c|c|c|c|c|c|c|c|}
        \hline
        Bound $b$ & Value & $b^2$ & Tol. & \# P/F & Av. \texttt{PASS} time [s] & Av. \texttt{FAIL} time [s] & Total time [s]\\\hline\hline
        $\Omega$ & $1.711\cdot 10^{-6}$ & $2^{-2q}$ & $10^{-14}$ & 4/8 & 218.1 & 982.8 & 8738.8 \\\hline
        $\epsilon$ & $6.8949\cdot 10^{-4}$ & $(4 \Omega/(\t \ss))^2$ & $10^{-9}$ & 11/4 & 73.8 & 75.2 & 1117.4 \\\hline
        $\delta$ & $1.383\cdot 10^{-3}$ & $2^{-q}$ & $10^{-9}$ & 10/4 & 84.7 & 90.8 & 1215.0 \\\hline
        $\Theta$ & $1.381\cdot 10^{-3}$ & $\delta^2$ & $10^{-9}$ & 3/8 & 12.1 & 28.3 & 266.6 \\\hline
        \end{tabular}
        \caption{Bounds obtained from the verification procedures for the three-mass-spring case study.} 
    \label{tab:bounds}
\end{table*}}

In our formal verification procedure we used the prototype tool of \cite{simic_tight_2022} for generating the bit-vector encoding, CBMC 5.4 \citep{clarke_CBMC_2004} for generating the SAT formula from the bit-vector encoding and MiniSat \citep{een_extensible_2004} to check for the satisfiability of the SAT formula.
All computations are performed on an Intel i5 processor running at $1.6$GHz.
We start by computing the value of $\Omega$ by selecting an initial value of $\Omega^2$ and then applying the bisection method on the verification of Algorithm~\ref{alg:verify:Omega} with an exit tolerance of $10^{-14}$, i.e., until the difference between the largest and smallest values of $\Omega^2$ resulting in a \texttt{FAIL} and a \texttt{PASS}, respectively, is smaller than $10^{-14}$.
Table~\ref{tab:bounds} shows the value of $\Omega$ obtained from this procedure, along with the selected exit tolerance, initial guess of $\Omega^2$, number of tests resulting in a \texttt{PASS}, number of tests resulting in a \texttt{FAIL}, average computation times of calls resulting in a \texttt{PASS} or \texttt{FAIL}, and total computation time of the bisection method. 
We take $\fp{\epsilon} = 2^{-q}$, which is the smallest value it can take, and then find the largest value of $\epsilon$ satisfying Assumption~\ref{alg:verify:ass:eps} by applying the bisection method on Algorithm~\ref{alg:verify:ass:eps}.
The results are presented in Table~\ref{tab:bounds}, where we note that the value of $\epsilon$ satisfies $\epsilon > 4 \Omega/(\t \ss) = 9.6195\cdot 10^{-5}$.
Therefore, we can use the exit condition $\distfp = 0$.
Finally, we obtain the bounds $\delta$ and $\Theta$ following the same bisection procedure used to compute $\Omega$.
The results are also presented in Table~\ref{tab:bounds}.
We take $\omega = \Omega$, as stated in Remark~\ref{rem:compute:omega}.

Plugging the results into Corollary~\ref{corr:certification}.\ref{corr:certification:dist}, we obtain the following: $k_\tmax = 250$, if $\distfp = 0$ then $\| \fpx^{k+1} - x^* \| \leq 0.0389$ and $f(\fpx^{k+1}) - f^* \leq 2.7165\cdot 10^{-4}$, are the best suboptimality bounds that can be guaranteed, since the ones from Corollary~\ref{corr:certification}.\ref{corr:certification:kmax} are smaller.
The value of $k_\tmax$ required to obtain the same $\epsilon/2$-suboptimality under exact arithmetic is $217$, c.f., \eqref{eq:k:exact}.

\vspace*{-1em}
\section{Conclusions} \label{sec:conclusions}

This article has presented a procedure for certifying the implementation of the PGM under fixed-point arithmetic when applied to strongly-convex box-constrained QP problems.
We have proven that the PGM maintains a linear convergence guarantee when sufficiently far away from the optimal solution, indicated by the choice of $\epsilon$, whose value can be reduced up to a maximum bound given by the fixed-point error-bound.
We have then presented a procedure based on recent formal verification tools to obtain a arbitrarily tight values of this error-bound and the other bounds that characterize the suboptimality of the output of the PGM.
Finally, we have shown that the computation times of the proposed verification procedures are tractable for a non-trivial MPC example.

\begin{appendix}
\gdef\thesection{\Alph{section}}
\makeatletter
\renewcommand\@seccntformat[1]{Appendix \csname the#1\endcsname.\hspace{0.5em}}
\makeatother

\section{Proofs and auxiliary lemmas} \label{app:proofs}

We start by providing two lemmas whose results are used in the proofs of Theorems~\ref{theo:convergence:fp} and \ref{theo:convergence:exit}.

\begin{lemma} \label{lemma:PGMfp:subgrad}
Consider Algorithm \ref{alg:PGMfp}.
For any $\alpha \in \Rpp$,
\begin{equation*}
\alpha \left( \fpx^{k} - \fpx^{k+1} - \gfp[\fpx^k] \right) \in \partial \Ind_{\cX}(\fpx^{k+1}), \; \forall k \geq 0.
\end{equation*}
\end{lemma}

\begin{proof}
    The claim follows from $\fpx^{k+1}$ being the Euclidean projection of $\fpx^k - \gfp[\fpx^k]$ onto $\cX$ (see Corollary \ref{corr:exact:projection}) along with the optimality condition of the projection operator \cite[Prop. 5.4.7]{bertsekas_convex_2009} and the equivalence between the subdifferential of the indicator function of a non-empty convex set and its normal cone \cite[Example 5.4.1]{bertsekas_convex_2009}.
\end{proof}

The following lemma particularizes \citep[Property~1]{alamo_restart_2019} to the fixed-point PGM paradigm.

\begin{lemma} \label{lemma:PGMfp:distf}
    Consider Algorithm \ref{alg:PGMfp} and let $v^k \in \B{\Omega}$ be the vectors that satisfy $\t \nabla f(\fpx^k) = \gfp[\fpx^k] + v^k$ for every $k \geq 0$.
    Denote $\stepk \doteq \fpx^{k} - \fpx^{k+1}$.
    Then, $\forall y \in \R^n$,
    {\def\theequation{\textit{\roman{equation}}}
        \begin{align}
            &f(\fpx^{k+1}) - f(y) \leq \ti \sp{\stepk + v^k}{\fpx^{k+1} - y} + \frac{1}{2 \t} \| \stepk \|^2 \label{eq:PGMfp:distf:1} \\
                                 &= \ti \sp{\stepk}{\fpx^{k} - y} - \frac{1}{2 \t} \| \stepk \|^2 + \ti \sp{v^k}{\fpx^{k+1} - y} \label{eq:PGMfp:distf:2} \\
                                 &= \frac{1}{2 \t} \| \fpx^k - y \|^2 {-} \frac{1}{2 \t} \| \fpx^{k+1} - y \|^2 + \ti \sp{v^k}{\fpx^{k+1} - y}. \label{eq:PGMfp:distf:3}
        \end{align}
    }
\end{lemma}

\begin{proof}
From Lemma~\ref{lemma:PGMfp:subgrad} we have that $\ti (\fpx^k - \fpx^{k+1} - \gfp[\fpx^k]) \in \partial \Ind_\cX(\fpx^{k+1})$.
Therefore, from the definition of the subdifferential \citep[\S 2.3]{parikh_neal_proximal_2013}, we have that
\begin{equation*}
    \Ind_\cX(y) \geq \Ind_\cX(\fpx^{k+1}) + \ti \sp{\fpx^k - \fpx^{k+1} - \gfp[\fpx^k]}{y - \fpx^{k+1}},
\end{equation*}
where taking $y \in \cX$ and recalling that $\fpx^{k+1} \in \cX$ (see Corollary~\ref{corr:exact:projection}), leads to
\begin{equation} \label{eq:proof:f_opt:1}
   0 \geq \ti \sp{\fpx^k - \fpx^{k+1} - \gfp[\fpx^k]}{y - \fpx^{k+1}}.
\end{equation}
From the convexity of $f$ we have that
\begin{equation} \label{eq:proof:f_opt:2}
    f(y) \geq f(\fpx^k) + \sp{\nabla f(\fpx^k)}{y - \fpx^k}.
\end{equation}
Additionally, from the $\L$-smoothness of $f$ and since $\t \leq \L^{-1}$, we have that
\begin{equation} \label{eq:proof:f_opt:3}
    f(\fpx^k) \geq f(\fpx^{k+1}) - \sp{\nabla f(\fpx^k)}{\fpx^{k+1} - \fpx^k} - \frac{1}{2 \t} \| \stepk \|^2.
\end{equation}
Claim~\eqref{eq:PGMfp:distf:1} follows from adding \eqref{eq:proof:f_opt:1}, \eqref{eq:proof:f_opt:2} and \eqref{eq:proof:f_opt:3} along with the definition of $v^k$.
Claims~\eqref{eq:PGMfp:distf:2} and \eqref{eq:PGMfp:distf:3} then follow from simple algebraic manipulations; c.f. Property 1.(i) in \cite{alamo_restart_2019}.\qedhere
\end{proof}

We now present the proof of Theorems~\ref{theo:convergence:fp}, which closely follows the proofs of \cite[Theorem 10.16 and Theorem 10.29]{beck_first-order_2017}, although various modifications have to be made to extend the results to the fixed-point arithmetic paradigm.

\begin{proof}[Proof of Theorem \ref{theo:convergence:fp}]
Let $\stepk \doteq \fpx^{k+1} - \fpx^k$ and $\dopt[k] \doteq \fpx^k - x^*$.
Consider the function $\psi : \cX \rightarrow \R$ given by
\begin{equation*}
    \psi(y) = f(\fpx^k) + \sp{\g{\fpx^k}}{y - \fpx^k} + \Ind_{\cX}(y) + \frac{1}{2 \t} \| y - \fpx^k \|^2.
\end{equation*}
Since $\psi$ is an $\ti$-strongly convex function, it follows from \cite[Theorem 5.24]{beck_first-order_2017} that
\begin{equation} \label{eq_proof_convergence_fp_psi_initial}
    \psi(y) - \psi(\fpx^{k+1}) \geq \sp{\mu}{y - \fpx^{k+1}} + \frac{1}{2\t} \| y - \fpx^{k+1} \|^2,
\end{equation}
$\forall \fpx^{k+1} \in \cX\pq$, $\forall \mu \in \partial \psi(\fpx^{k+1})$.
Since $\t \leq \L^{-1}$, we have that $f(\cdot)$ satisfies the well-known descent lemma \cite[Lemma 5.7]{beck_first-order_2017}
\begin{equation*}
f(\fpx^{k+1}) \leq f(\fpx^k) + \sp{\g{\fpx^k}}{\fpx^{k+1} - \fpx^k} + \frac{1}{2\t} \| \fpx^{k+1} - \fpx^k \|^2,
\end{equation*}
$\forall \fpx^{k+1} \in \cX\pq$, $\forall y \in \cX$, which along with the definition of $\psi$ and noting that $\Ind_{\cX}(\fpx^{k+1}) = 0$ by virtue of Corollary \ref{corr:exact:projection}, leads~to
\begin{equation*}
    \psi(\fpx^{k+1}) = f(\fpx^k) + \sp{\nabla f(\fpx^k)}{\stepk} + \frac{L}{2} \| \stepk \|^2 \geq f(\fpx^{k+1}), 
\end{equation*}
$\forall \fpx^{k+1} \in \cX\pq$.
Thus, we can rewrite \eqref{eq_proof_convergence_fp_psi_initial} as
\begin{equation} \label{eq:proof:fp:psi:ineq}
    \psi(y) - f(\fpx^{k+1}) \geq \sp{\mu}{y - \fpx^{k+1}} + \frac{1}{2\t} \| y - \fpx^{k+1} \|^2,
\end{equation}
$\forall \fpx^{k+1} \in \cX\pq, \forall \mu \in \partial \psi(\fpx^{k+1})$.
The subdifferential of $\psi$ evaluated at $\fpx^{k+1}$ is given by
\begin{equation} \label{eq:subdif:psi}
\partial \psi(\fpx^{k+1}) = \g{\fpx^k} + \ti (\fpx^{k+1} - \fpx^{k}) + \partial \Ind_{\cX}(\fpx^{k+1}).
\end{equation}
From Definition~\ref{def:Omega} we have that for each $k \geq 0$ there exists a vector $v^k \in \B[n]{\Omega}$ satisfying $\t \g{\fpx^k} = \gfp + v^k$.
Therefore, we can rewrite \eqref{eq:subdif:psi} as
\begin{equation*}
\partial \psi(\fpx^{k+1}) = \ti (v^k + \gfp + \fpx^{k+1} - \fpx^{k} ) + \partial \Ind_{\cX}(\fpx^{k+1}).
\end{equation*}
From Lemma~\ref{lemma:PGMfp:subgrad} we have $0 \in \ti ( \gfp + \fpx^{k+1} - \fpx^{k}) + \partial \Ind_{\cX}(\fpx^{k+1})$, thus $\ti v^k \in \partial \psi(\fpx^{k+1})$.
This allows us to rewrite \eqref{eq:proof:fp:psi:ineq} as
\begin{equation*}
    \psi(y) - f(\fpx^{k+1}) \geq \sp{L v^k}{y - \fpx^{k+1}} + \frac{1}{2\t} \| y - \fpx^{k+1} \|^2,
\end{equation*}
$\forall \fpx^{k+1} \in \cX\pq$, for some $v^k \in \B[n]{\Omega}$.
Undoing the expression of $\psi(y)$ and particularizing to $y = x^* \in \cX$ leads to
\begin{align*}
    f(\fpx^*) - f(\fpx^{k+1}) \geq& \frac{1}{2\t} \| \dopt[k+1] \|^2 - \frac{1}{2\t} \| \dopt[k] \|^2 + f(\fpx^*) - f(\fpx^k) \\
    &- \sp{\nabla f(\fpx^k)}{\fpx^* - \fpx^k} + \ti \sp{v^k}{\fpx^* - \fpx^{k+1}}.
\end{align*}
Since $f$ is a $\ss$-strongly convex function, we have that \cite[Theorem 5.24.\textit{(ii)}]{beck_first-order_2017}
\begin{equation*}
    f(\fpx^*) - f(\fpx^k) - \sp{\nabla f(\fpx^k)}{\fpx^* - \fpx^k} \geq \frac{\ss}{2} \| \fpx^* - \fpx^k \|^2.
\end{equation*}
Thus,
\begin{equation} \label{eq:proof:fp:f:bound}
    f(\fpx^*) - f(\fpx^{k+1}) \geq \frac{1}{2\t} \| \dopt[k+1] \|^2 - \frac{\ti - \ss}{2} \| \dopt[k] \|^2 - \ti \sp{v^k}{\dopt[k+1]}.
\end{equation}
By definition of $\fpx^*$, we have that $f(\fpx^*) - f(\fpx^{k+1}) \leq 0$.
Therefore, the right hand side of \eqref{eq:proof:fp:f:bound} must also be less or equal to $0$, which leads to
\begin{align*}
    \frac{1}{2} \| \fpx^{k+1} - \fpx^* \|^2 &\leq \frac{1 - \t\ss}{2} \|\fpx^k - \fpx^* \|^2 + \sp{v^k}{\fpx^{k+1} - \fpx^*} \\
                                            &\becauseof[\leq]{(*)} \frac{1 - \t\ss}{2} \|\fpx^k - \fpx^* \|^2 + \| v^k \| \cdot \| \fpx^{k+1} - \fpx^* \| \\
                                            &\becauseof[\leq]{(**)} \frac{1 - \t\ss}{2} \|\fpx^k - \fpx^* \|^2 + 2 \Omega \epsilon^{-1} \| \fpx^{k+1} - \fpx^* \|^2,
\end{align*}
where in $(*)$ we are making use of the Cauchy-Schwarz inequality and in $(**)$ of the fact that $\| v^k \| \leq \Omega$ and $\| \fpx^{k+1} - \fpx^* \| > \fracg{\epsilon}{2}$.
Since by construction $\t \ss \leq 1$, the assumption $4 \Omega \epsilon^{-1} < \t \ss$ implies $2 \Omega \epsilon^{-1} < \fracg{1}{2}$.
Thus, we derive
\begin{align*}
    \| \fpx^{k+1} - \fpx^* \|^2 &\leq \left( \fracg{1 - \t \ss }{1 - 4 \Omega \epsilon^{-1}} \right) \| \fpx^k - \fpx^* \|^2,
\end{align*}
which leads to claim~\ref{theo:convergence:fp:xk} if applied recursively, where we note that the assumption $4 \Omega \epsilon^{-1} < \t \ss$ guarantees that the sequence is convergent.

We now prove claim~\ref{theo:convergence:fp:f} by rearranging \eqref{eq:proof:fp:f:bound} and proceeding as follows:
\begin{align*}
    f(\fpx^{k+1}) &{-} f(\fpx^*) \leq \frac{\ti - \ss}{2} \| \dopt[k] \|^2 - \frac{1}{2\t} \| \dopt[k+1] \|^2 + \ti\sp{v^k}{\dopt[k+1]} \\
                                &\leq \frac{\ti - \ss}{2} \| \dopt[k] \|^2 - \frac{1}{2\t} \| \dopt[k+1] \|^2 + 2 \ti \Omega \epsilon^{-1} \|\dopt[k+1]\|^2 \\
                                &\leq \frac{\ti - \ss}{2} \| \dopt[k] \|^2 + \ti \left( 2 \Omega \epsilon^{-1} - \frac{1}{2} \right) \| \dopt[k+1] \|^2 \\
                                &\becauseof[\leq]{(*)} \frac{ 1 - \t\ss }{2\t}  \| \dopt[k] \|^2 = \frac{1 - 4 \Omega \epsilon^{-1}}{2\t} \left( \frac{1 - \t\ss}{1 - 4 \Omega \epsilon^{-1}} \right) \| \dopt[k] \|^2 \\
                              &\becauseof[\leq]{(**)} \frac{1 - 4 \Omega \epsilon^{-1}}{2\t} \left( \frac{1 - \t\ss }{1 - 4 \Omega \epsilon^{-1}} \right)^{k+1} \| \fpx^0 - \fpx^* \|^2,
\end{align*}
$\forall \fpx^{k+1} \in \cX\pq$, $\forall k \geq 0$,
where $(*)$ holds since $2 \Omega \epsilon^{-1} < \fracg{1}{2}$ and $(**)$ follows from claim~\ref{theo:convergence:fp:xk}.\qedhere
\end{proof}

\end{appendix}

\bibliographystyle{elsarticle-harv}
\bibliography{CertQP_PGM_bib}

\end{document}